\documentclass[12pt,reqno]{article}

\usepackage[pdftex]{pict2e}
\usepackage{cite,euscript,verbatim,afterpage,bm,paralist,stmaryrd}
\usepackage{lineno}
\usepackage{hyperref}

\usepackage{amsmath}
\usepackage{amssymb}
\usepackage{amsthm,mathtools}
\usepackage{amscd}
\usepackage{amsfonts}
\usepackage{stfloats}
\usepackage{float}
\usepackage{graphicx}
\usepackage{color}
\usepackage{graphicx,url,newverbs}
\usepackage{caption}
\usepackage[ruled,vlined]{algorithm2e}
\usepackage{algorithmic}
\usepackage{verbatim}
\usepackage{enumitem}

\newcommand{\fig}[1]{\mbox{Figure{~#1}}}

\numberwithin{equation}{section}

\newcommand{\R}{{\mathbb{R}}}

\newcommand{\bx}{{\mathbf{x}}}
\newcommand{\bmu}{{\boldsymbol{\mu}}}
\newcommand{\blambda}{{\boldsymbol{\lambda}}}
\newcommand{\bu}{\mathbf{u}}

\newcommand{\bpsi}{\boldsymbol{\psi}}
\newcommand{\bpoly}{\mathbf{p}}
\newcommand{\half}{\frac{1}{2}}

\newcommand{\icol}[1]{
	\left(\begin{smallmatrix}#1\end{smallmatrix}\right)%
}

\newcommand{\ba}{\begin{array}}
\newcommand{\ea}{\end{array}}
\newcommand{\be}{\begin{equation}}
\newcommand{\ee}{\end{equation}}
\newcommand{\bd}{\begin{displaymath}}
	\newcommand{\ed}{\end{displaymath}}

\begin{document}

\allowdisplaybreaks

\title{
Integration factor combined with level set method for reaction-diffusion systems with free boundary in high spatial dimensions}

\author{
Shuang Liu
\thanks{Department of Mathematics, University of California, San Diego, 
9500 Gilman Drive, La Jolla, California 92093-0112, United States.
Email: shl083@ucsd.edu
}
\and
Xinfeng Liu
\thanks{
Corresponding Author, Department of Mathematics, University of South Carolina, Columbia, South Carolina 29208, United States. 
Email: xfliu@math.sc.edu. This work is partially supported by NSF DMS1853365.
}
}


\maketitle

\begin{abstract}
\noindent
For reaction-diffusion equations in irregular domain with moving boundaries, the numerical stability constraints from the reaction and diffusion terms often require very restricted time step size, while complex geometries may lead to difficulties in accuracy when discretizing the high-order derivatives on grid points near the boundary. It is very challenging to design numerical methods that can efficiently and accurately handle both difficulties.  Applying an implicit scheme may be able to remove the stability constraints on the time step, however, it usually requires solving a large global system of nonlinear equations for each time step, and the computational cost could be significant. Integration factor (IF) or exponential differencing time (ETD) methods 
are one of the popular methods for temporal partial differential equations (PDEs) among many other methods. In our paper, we couple ETD methods with an embedded boundary method to solve a system of  reaction-diffusion equations with complex geometries. In particular, we rewrite all ETD schemes into a linear combination of specific $\phi$-functions and apply one start-of-the-art algorithm  to compute the matrix-vector multiplications, which offers significant computational advantages with adaptive Krylov subspaces. In addition, we extend this method by incorporating  the level set method to solve the free boundary problem. The  accuracy, stability, and efficiency of the developed method are demonstrated by numerical examples.

\bigskip

\noindent
{\bf Keywords:}
Reaction Diffusion Equations, Free Boundary, Integrating Factor Method, Level Set Method

\end{abstract}

\section{Introduction}

The systems of reaction-diffusion equations coupled with moving boundaries defined by Stefan condition have been widely used to describe the dynamics of the spreading population.
A moving boundary problem is characterized by the fact that the boundary of the domain is not known in advance but it has to be determined as a part of the solution. These problems are often called Stefan problems due to the Stefan condition that links the behavior of the boundary with the unknown solution \cite{piqueras2017front,rubinshteuin1971stefan,tayler1985free}.

The Stefan condition was first introduced with a moving boundary of parabolic type to describe the spreading of species population as introduced in \cite{du2010spreading}, the reaction-diffusion system for the density of population of the invasive species $u(\mathbf{x},t)$ depending on time $t$ and spatial variable $\mathbf{x}$. In this paper, we consider solving the following system of reaction-diffusion equations coupled with free boundaries,  
\begin{equation}\label{eq:free}
	\left \{
	\begin{aligned}
		& 	\frac{\partial u}{\partial t}  =\nabla\cdot(\beta \nabla {u})+f({{u}},t)&\mathbf{x}\in \Omega(t),\\
		&{u}({\mathbf{x}},t)=0&\mathbf{x}\in \Omega^{c}(t),\\
		&\vec{v}(\mathbf{x},t)=-\mu\nabla u\cdot \mathbf{n},&\mathbf{x}\in \partial\Omega(t).
	\end{aligned}
	\right.
\end{equation}
where $\mu>0$, $\partial\Omega(t) $ is the moving boundary of the  evolution of the  domain $\Omega(t)$, which represents the spreading front of the species $u({\mathbf{x}},t)$. 
Here the evolution of the moving domain $\Omega(t)\subset \mathbb R^N$, or rather its boundary $\partial\Omega(t)$ is determined by the one phase Stefan condition
which, in the case $\partial\Omega(t)$ is a $C^1$ manifold in $\mathbb R^N$, can be described as follows:
\begin{description}
\item \hspace{0.8cm}Any point ${\bf x} \in \partial\Omega(t)$ moves with velocity $\mu |\nabla_{\bf x} u({\bf x},t)|{\bf n}({\bf x})$, where ${\bf n}({\bf x})$
\\
 \hspace{1.2cm}  is the unit outward normal of $\Omega(t)$ at ${\bf x}$, and $\mu$ is a given positive constant.
\end{description}

The moving boundary is generally called the ``free boundary'', which has been extensively studied theoretically \cite{caffarelli2005geometric} and numerically \cite{cao1989numerical,chen2009numerical,chen1997simple, gibou2005fourth,OshFed02,osher2001level,peng1999pde} and the references therein. Other theoretical studies of related free boundary problems can be found in \cite{bunting2012spreading} and the references therein.
When solving such a system numerically, difficulties arise from the stiffness along with moving boundaries. 
First of all, it is always extremely difficult to handle points near the boundary.
To overcome this, various numerical techniques have been proposed for providing the pros and cons of different choices for defining the ghost values to avoid the small cell stiffness, while those numerical treatments focus on introducing a small positive number as  the threshold of the distance between the interior points and the boundary points \cite{gibou2005fourth,gibou2002second,jomaa2005embedded,  liu2000boundary}.  To some extent, these techniques can remove the large errors that could occur from dividing by small numbers to get second-order accurate solutions, however, remedies are required to keep the numerical accuracy of the gradients by not only  proposing higher order extrapolating for defining the ghost points, but also  combining higher order interpolation for locating the interface \cite{ng2009guidelines}. 

To overcome this difficulty, we adopt an embedded boundary method to solve a variable coefficient Poisson equation in an irregular domain with Dirichlet boundary conditions. Numerical solutions to the Poisson equation in irregular domains have been considered by many approaches, including finite difference \cite{buzbee1971direct, gibou2002second,jomaa2005embedded, lai2001note, min2006supra,perrone1975general,yoon2015analyses}, finite volume \cite{ johansen1998cartesian,oevermann2006cartesian,zapata2018gpu}, and finite element \cite{apel1996graded, braess1981contraction, louis1979acceleration, sudicky1989laplace} using various meshing techniques. Among them, the embedded boundary method has a number of advantages, which includes simplifying the grid generation process for complicated geometries, enabling fast computation approach in parallel, and shifting the complexity of dealing with complex geometries to the discretization approach. More importantly, the embedded boundary method is an excellent candidate with extension to the moving boundary problems, as it generates the mesh using a background regular mesh by taking special care of cut-cells where the geometry intersects the grid.

The placement of the \emph{ghost point} is the subtle yet important distinction from a wide range of methods \cite{Kreiss_Petersson_2006,kreiss:1940,kreiss:1292}. In contrast to the Ghost Fluid Method introduced by \cite{gibou2002second}, where \emph{ghost points} are  placed outside the computational domain, here we plan to use  \emph{interior ghost points} instead to ease the small cell stiffness when the interface is very close to the grid points in the irregular domain. The proposed embedded boundary method results in a symmetric positive definite  discretization matrix, thus we can use a wild number of fast linear solvers. For instance, algebraic multigrid  with both ``V"-cycle and  ``W"-cycle can be applied as preconditioners to further speed up calculations.

On the other hand, extremely small time steps are required due to the stiffness of the system. When the explicit schemes are applied to solve such a system, due to stability constraints, an extremely small time step should be used and it might take a long time to finish one single simulation. However, while applying an implicit scheme~\cite{BurBut79, HaiWan99, Low04} may be able to remove the stability constraints on the time step $\Delta t$, it usually requires solving a large global system of nonlinear equations for each time step, and the computational cost could be significant.

To remove the stability constrains on the size of time steps,  we employ  exponential time differencing (ETD) methods in which the diffusion term is discretized by the embedded boundary method. As is known, the ETD schemes exhibit very nice stability properties, which allow for the relatively large time step size~\cite{BeyKei98, DuJu21, DuZhu04a, DuZhu05, HouLow94, jiang2021unconditionally, JouLeo97, KasTre05, LeoLow00, li2021unconditionally}. In addition, by rewriting all ETD schemes as a linear combination of $\phi$-functions, we combine a state-of-the-art algorithm: $phipm_{-}simul_{-}iom2$  \cite{luan2019further} to evaluate the linear combination of matrix-vector  multiplications, which offers a significant computational advantages by adopting adaptive Krylov subspaces.

The rest of the paper is organized as follows. In \S 2, a second-order embedded boundary method is presented to discretize the diffusion term in irregular domains. In \S 3, we briefly describe the explicit ETD schemes along with Runge-Kutta type of ETD schemes by rewriting all ETD schemes as a linear combination of $\phi$-functions. Furthermore, we introduce a state-of-the-art algorithm for computing linear combinations of matrix function $\phi_k(A)$ on vectors $v_k$.  In \S 4, various numerical examples have been performed to demonstrate the accuracy, efficiency, and stability of the developed algorithms, and such schemes have been incorporated with level set method  to solve the  free boundary problem as depicted in~(\ref{eq:free}). Finally, in \S 5, we draw a brief conclusion and further discuss several possible extensions for future studies.

\section{A Cartesian Grid Embedded Boundary Method }\label{sec:embedded}
In this section, following \cite{peng2022universal}, we briefly  introduce a Cartesian grid embedded boundary method to develop a second-order symmetric positive definite discretization of a static Poisson equation with Dirichlet boundary conditions in irregular domains.

In order to illustrate the approach, we first consider a stationary Poisson equation with variable coefficients in a two-dimensional domain $\Omega$
\begin{equation} \label{eq:peq} 
\nabla\cdot(\beta(x,y)\nabla u) = f(x,y), \quad (x,y)\in \Omega,
\end{equation}
with Dirichlet boundary conditions
\begin{equation}
\label{eq:peqbcD} u(x,y) = u_{\mathcal{D}}(x,y),\quad (x,y)\in \partial \Omega. 
\end{equation}
on the interface. Without loss of generality, we assume that the irregular domain $\Omega$ are contained inside a rectangular domain
$[x_L,x_R] \times [y_L,y_R]$, covered by a uniform grid with a grid function denoted by $u_{ij}=u(x_i,y_j)$, where
\begin{equation}
\bx_{ij}=(x_i,y_j) = (x_L+(i-1)h,y_L+(j-1)h),\ \ i=1,\ldots,n_x,\,\, j=1,\ldots,n_y.
\end{equation}

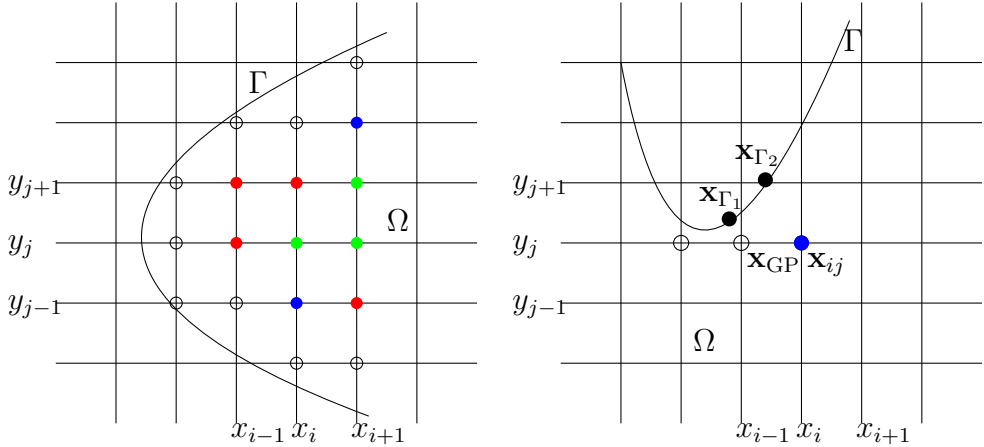
\begin{figure}[htb]
\begin{center}
\setlength{\unitlength}{0.8cm} 
\begin{picture}(7,7) 
\multiput(0,1)(0,1){2}{\line(1,0){7}}
\multiput(0,4)(0,1){3}{\line(1,0){7}}
\multiput(1,0)(1,0){6}{\line(0,1){7}}
\qbezier(5.5,6.5)(-2.5,3)(5.2,0.12)
\put(2,3){\circle{0.2}}
\put(2,4){\circle{0.2}}
\put(3,5){\circle{0.2}}
\put(2,2){\circle{0.2}}
\put(3,2){\circle{0.2}}
\put(4,5){\circle{0.2}}
\put(4,1){\circle{0.2}}
\put(5,6){\circle{0.2}}
\put(3.2,5.5){$\Gamma$}
\put(5.5,3.2){$\Omega$}

\put(0.,3){\line(1,0){7}}
\put(3.91,-0.23){$x_i$}
\put(4.91,-0.23){$x_{i+1}$}
\put(2.91,-0.23){$x_{i-1}$}
\put(-0.8,2.9){$y_j$}
\put(-0.8,1.9){$y_{j-1}$}
\put(-0.8,3.9){$y_{j+1}$}
\put(5,1){\circle{0.2}}
\color{red}{\put(3,3){\circle*{0.2}} } 
\put(4,4){\circle*{0.2}}
\put(3,4){\circle*{0.2}}
\put(5,2){\circle*{0.2}}
\color{green}
\put(5,4){\circle*{0.2}}
\put(5,3){\circle*{0.2}}
\put(4,3){\circle*{0.2}}
\color{blue}
\put(5,5){\circle*{0.2}}
\put(4,2){\circle*{0.2}}
\end{picture}
\hspace{0.7cm}
\setlength{\unitlength}{0.8cm} 
\begin{picture}(7,7) 
\multiput(0,1)(0,1){2}{\line(1,0){7}}
\multiput(0,4)(0,1){3}{\line(1,0){7}}
\multiput(1,0)(1,0){6}{\line(0,1){7}}
\qbezier(4.8,6.7)(2.1,0.1)(1,6)
\put(4,3){\circle*{0.25}}
\put(3,3){\circle{0.25}}
\put(2,3){\circle{0.25}}
\put(3.4,4.05){\circle*{0.25}}
\put(2.8,3.4){\circle*{0.25}}

\put(2.9,4.4){$\bx_{\Gamma_2}$} 
\put(2.3,3.7){$\bx_{\Gamma_1}$}

\put(4.7,6.2){$\Gamma$}

\put(2.2,1.2){$\Omega$}
\put(3.91,-0.23){$x_i$}
\put(4.91,-0.23){$x_{i+1}$}
\put(2.91,-0.23){$x_{i-1}$}
\put(-0.8,2.9){$y_j$}
\put(-0.8,1.9){$y_{j-1}$}
\put(-0.8,3.9){$y_{j+1}$}

\put(0,3){\line(1,0){7}}
\put(4.1,2.6){$\bx_{ij}$} 
\put(3.1,2.6){$\bx_{\textrm{GP}}$} 

\color{blue}\put(4,3){\circle*{0.25}}

\end{picture}
\caption{Left: Discretization in an irregular domain  with circles representing {\em interior ghost points} and colored dots representing {\em computational points}. Right: Illustration to construct a RBF interpolation to obtain $u_{\textrm{GP}}$.\label{fig:LinebylineDis}}
\end{center}
\vspace{-10mm}
\end{figure}

We denote grid points inside $\Omega$ on the fringe of the
computational domain  as \emph{interior ghost points}, and all other interior grid points are denoted as \emph{computational points}. An \emph{interior ghost point} $(x_i,y_j)$ satisfies the condition that $(x_i,y_j)$ is inside $\Omega$, but at least one of its four nearest neighbors is outside. While, a \emph{computational point} $(x_i,y_j)$ satisfies the condition that $(x_i,y_j)$ and all its four nearest neighbors are all inside $\Omega$. For example, in \fig{\ref{fig:LinebylineDis}} (left), all colored grid points are \emph{computational points} while all grid points with a circle are \emph{interior ghost  points}.   Note that the Poisson system \eqref{eq:peq} will be solved only at  \emph{computational points} while not at \emph{interior ghost points}.

There are three types of {\em{computational points}} to discretize the Laplacian operator as follows:

{\bf Case 1.} If the {\em{computational point}} is in the  absence of {\em{interior ghost points}} as its neighbors (green grid points in \fig{\ref{fig:LinebylineDis}} (left)), Laplacian operator  is approximated by a standard central difference scheme. We take the point $(x_i,y_j)$ in \fig{\ref{fig:LinebylineDis}} (left) as an example, 
\begin{align}
\nabla\cdot(\beta(x_i,y_j)\nabla u(x_i,y_j))&\approx 
\frac{\beta_{i+\half,j}u_{i+1,j}-(\beta_{i+\half,j}+\beta_{i-\half,j})u_{i,j}+\beta_{i-\half,j}u_{i-1,j}}{h^2}\notag\\
&+\frac{\beta_{i,j+\half}u_{i,j+1}-(\beta_{i,j+\half}+\beta_{i,j-\half})u_{i,j}+\beta_{i,j-\half}u_{i,j-1}}{h^2}.
\label{eq:central_difference}
\end{align}

{\bf Case 2.} For one coordinate direction, the {\em{computational point}} neighbors an {\em{interior ghost point}} while the {\em{interior ghost point}} borders the interface, for example, the red and blue grid points in \fig{\ref{fig:LinebylineDis}} (left) for $y$ coordinate direction. The Lagrange polynomial interpolation with a line by line approach will be applied for this case.

For instance, we take $(x_i,y_{j+1})$ in $y$ direction as an example for illustration (See \fig{\ref{fig:LinebylineDis}} (left)). 
The intersection point of the grid line $x=x_i$ with the boundary $\Gamma$ between $(x_{i},y_{j+2})$ and $(x_{i},y_{j+3})$  is denoted by $(x_i,y_\Gamma)$, and the boundary value at  $(x_i,y_\Gamma)$ is given by  $u_{\Gamma}$. Here $y_\Gamma$ can be found by some root-finding algorithm, such as the secant method.

Next we introduce an interpolation polynomial such that  the value of the {\em{interior ghost point}} $u_{\textrm{GP}}=u_{i,j+2}$ can be estimated as
\begin{equation*}
\label{eq:DirichletLineByLine} u_{\textrm{GP}}=\mathcal{I}_{P}u(x_i,y_{j+2}) = u_\Gamma g_\Gamma(y_{j+2})+u_{i,j+1} g_1(y_{j+2})= u_\Gamma \frac{y_{j+2}-y_{j+1}}{y_\Gamma-y_{j+1}}-u_{i,j+1}\frac{y_{j+2}-y_\Gamma}{y_\Gamma-y_{j+1}}.
\end{equation*}

Substituting $u_{\textrm{GP}}$ into the central difference approximation for the Laplacian operator in $y$ direction at the point $(x_i,y_{j+1})$, we obtain

\begin{equation*}
\frac{-(\beta_{i,j+\frac{3}{2}}(1-g_1(y_\textrm{GP}))+\beta_{i,j+\half})u_{i,j+1}+\beta_{i,j+\half}u_{i,j}}{h^2}+\frac{\beta_{i,j+\frac{3}{2}}u_{\Gamma}g_{\Gamma}(y_\textrm{GP})}{h^2}.
\label{eq:central_difference2}
\end{equation*}

Noticing that
\begin{displaymath}
|g_{\Gamma}(y_{\textrm{GP}})|=\left|\frac{y_{\textrm{GP}}-y_{j+1}}{y_\Gamma-y_{j+1}}\right|=\left|\frac{y_{j+2}-y_{j+1}}{y_\Gamma-y_{j+1}}\right|\leq 1\;\text{and}\; |g_1(y_{\textrm{GP}})|=\left|\frac{y_{\textrm{GP}}-y_\Gamma}{y_\Gamma-y_{j+1}}\right|=\left|\frac{y_{j+2}-y_\Gamma}{y_\Gamma-y_{j+1}}\right|\leq 1.
\end{displaymath}

The resulting linear system is still diagonal dominant with correct sign, and the symmetric structure is also preserved with only diagonal elements modified.

{\bf\ Case 3.} For one coordinate direction, the {\em{computational point}} neighbors an {\em{interior ghost point}} while this {\em{interior ghost point}} does not  border the interface, for example, the blue grid points in \fig{\ref{fig:LinebylineDis}} (left) for $x$ coordinate direction. For this case, a Radial basis function (RBF) based interpolation \cite{franke1982scattered, hardy1971multiquadric} will be employed since applying the Lagrange polynomial interpolation directly would cause loss of accuracy.

Without loss of generality, we consider the case in a non-convex geometry as presented in \fig{\ref{fig:LinebylineDis}} (right) that will occur in later numerical testing examples.  We take $(x_i,y_j)$ in $x$ coordinate direction as an example for illustration.  Let $\bx_{ij}=(x_i,y_j)$ be a {\em{computational point}} and $\bx_{\textrm{GP}}=(x_{i-1},y_{j})$ be an  {\em{interior ghost  point}} neighboring $\bx_{ij}$ in $x$ direction. 

We choose $\bx_{\Gamma_1}$ and $\bx_{\Gamma_2}$  to be the corresponding closet points on the boundary to $\bx_{\textrm{GP}}$ and $\bx_{ij}$. 
We use the following combination of RBF and a linear polynomial tail for interpolation at $\bx_{\textrm{GP}}$:
\begin{align*}
\mathcal{I}_{\textrm{R}}u(\bx) = \lambda_{ij} \psi(||\bx-\bx_{ij}||)+\lambda_{\Gamma_1}\psi(||\bx-\bx_{\Gamma_1}||)+\lambda_{\Gamma_2}\psi(||\bx-\bx_{\Gamma_2}||)+\mu_1+\mu_2 x+\mu_3 y,  
\end{align*}
where $\bx=(x,y)$, $||\cdot||$ is the standard $l_2$ norm and $\psi(\cdot)$ is a radial basis function. The linear polynomial tail is required to maintain second order accuracy \cite{barnett2015robust}. The coefficient $\blambda=(\lambda_{ij},\lambda_{\Gamma_1},\lambda_{\Gamma_2})^T$ and $\bmu=(\mu_1,\mu_2,\mu_3)^T$ are determined by the linear system

\be
B\left(\begin{matrix}
\blambda \\
\bmu 
\end{matrix}\right) 
=\left(\begin{matrix}
A & \Pi^T \\
\Pi & 0 
\end{matrix}\right) 
\left(\begin{matrix}
\blambda \\
\bmu 
\end{matrix}\right) 
=\left(\begin{matrix}
\bu \\
{\mathbf{0}}
\end{matrix}\right).
\ee 

where $\bu=(u_{ij},u_{\Gamma_1},u_{\Gamma_2})$,
\bd
A = \left(\begin{matrix}
\psi(0) & \psi(||\bx_{ij}-\bx_{\Gamma_1}||) &  \psi(||\bx_{ij}-\bx_{\Gamma_2}||) \\
\psi(||\bx_{\Gamma_1}-\bx_{ij}||) & \psi(0) &  \psi(||\bx_{\Gamma_1}-\bx_{\Gamma_2}||) \\
\psi(||\bx_{\Gamma_2}-\bx_{ij}||) &  \psi(||\bx_{\Gamma_2}-\bx_{\Gamma_1}||) &  \psi(0) 
\end{matrix}\right) 
\quad
\text{and}
\quad
\Pi = \left(\begin{matrix}
1 & 1 & 1\\
x_{ij} & x_{\Gamma_1} & x_{\Gamma_2} \\ 
y_{ij} & y_{\Gamma_1} & y_{\Gamma_2}
\end{matrix}\right).
\ed
The  value at the {\em{boundary point}} $\bx_{\textrm{GP}}$ is assigned as
\begin{align}
u_{\textrm{GP}} =\mathcal{I}_{\textrm{RBF}}u(\bx_{\textrm{GP}}) = \left(\bpsi_{\textrm{GP}}^T,\bpoly_{\textrm{GP}}^T\right)
B^{-1}\left(\begin{matrix}
\bu \\
{\mathbf{0}}
\end{matrix}\right)
\label{eq:rbf_ubp}
\end{align} 
with $\bpsi_{\textrm{GP}}=\left(\psi(||\bx_{\textrm{GP}}-\bx_{ij}||), \psi(||\bx_{\textrm{GP}}-\bx_{\Gamma_1}||), \psi(||\bx_{\textrm{GP}}-\bx_{\Gamma_2}||)\right)^T$ and $\bpoly_{\textrm{GP}}=(1,x_{\textrm{GP}},y_{\textrm{GP}})^T$. As the terms in the right hand side in \eqref{eq:rbf_ubp} appear only in the diagonal coefficient, the symmetry of the discrete matrix will not be broken. For this case, the Laplacian operator in $x$ direction at the point $(x_i,y_j)$ can be approximated by 
\begin{equation*}
\frac{\beta_{i+\half,j}u_{i+1,j}-(\beta_{i+\half,j}+\beta_{i-\half,j})u_{i,j}+\beta_{i-\half,j}u_{\textrm{GP}}}{h^2}.
\label{eq:central-difference3}
\end{equation*}

\section{Exponential Time Differencing Schemes}
In this section, we briefly discuss both the explicit ETD schemes and  Runge-Kutta type of ETD schemes with arbitrary order accuracy. For illustration, here we consider a reaction-diffusion system with certain boundary conditions,
\begin{equation}\label{Eq01}
	\mathbf{u}_t=\nabla\cdot(\beta(x,y)\nabla \mathbf{u})+\mathbf{f}(\mathbf{u},t),
\end{equation}
where $\mathbf{u}=\mathbf{u}(\mathbf{x},t)\in R^m$, $\beta(x,y)$ is the diffusion coefficient, and $\mathbf{f}(\mathbf{u},t)$ represents nonlinear reaction terms.

Applying the spatial discretization with the embedded boundary method as previously illustrated in \S \ref{sec:embedded} for the Poisson equation,  we reduce the equation (\ref{Eq01}) to a system of ODEs: 
\begin{equation}\label{Eq02}
	U_t=CU+F(U(t),t),
\end{equation}
where  $U=U(t)$ is the spatially discretized form of $\mathbf{u}$,  and $C$ is a constant matrix representing the finite difference approximation of  the diffusion. After multiplying the equation (\ref{Eq02}) by the integrating factor $e^{-Ct}$, we integrate the equation over one time step from $t_n$ to $t_{n+1}\equiv t_n+\triangle t$ to obtain 
\begin{equation}\label{Eq03}
	U(t_{n+1})=U(t_{n})e^{C\triangle t}+\int_0^{\triangle t}e^{(\triangle t-\tau)C}F(U(t_n+\tau), t_n+\tau)d\tau.
\end{equation}
While this formula is exact, the essence of the ETD methods is to derive numerical approximations to the integral in this expression. 

\subsection{Explicit exponential time differencing }

For the derivation of ETD schemes, the integrand is approximated first through interpolation polynomials of the function $F(U(t_n+\tau), t_n+\tau)$ with $e^{-C\tau}$ unchanged. With the Lagrange interpolation being applied to 
approximate $F(U(t_n+\tau), t_n+\tau)$, a direct integration of the interpolation polynomial with the coefficient term $e^{-C\tau}$ yields the ETD method. If all interpolation points used for the integrand are with $\tau \leq 0$, 
the resulted temporal scheme is explicit. Otherwise, the scheme becomes implicit when the interpolation points contain the solution at $t_{n+1}$. 


Assuming that $F(U(t),t)$ is constant such that $F(U(t),t) =F_n\equiv F(U_n,t_n)$ over the interval $t_n\leq t\leq t_{n+1}$, we introduce a scalar function
\begin{equation*}
	\phi_1(z)=\frac{e^z-1}{z},
\end{equation*}
and the first order ETD scheme (ETD1) is given by
\begin{align*}\label{Eq04}
	U_{n+1}&=U_{n}e^{C\triangle t}+\triangle t \phi_1(\triangle t C)F_n\\
	&=U_n+\triangle t \phi_1(\triangle t C)(F_n+CU_n),
\end{align*}
where $U_{n+1}$ is the numerical approximation to $U(t_{n+1})$ and $U(t_n)$ as $U_n$. Here
\begin{equation*}
	\phi_1(\triangle t C)=\frac{1}{\triangle t}\int_{0}^{\triangle t}e^{(\triangle t-\tau)C}d\tau=\int_{0}^{1}e^{(1-\lambda)\triangle t C}d\lambda
\end{equation*}
is extended to the matrix form from the scalar function $\phi_1$.

For the second-order approximation
\begin{equation*}\label{Eq05}
	F(U(t_n+\tau), t_n+\tau)=F_n+\tau\frac{F_n-F_{n-1}}{\triangle t}+O(\triangle t)^2,
\end{equation*}
the second-order ETD scheme (ETD2) can be achieved by 
\begin{align*}\label{Eq06}
	U_{n+1}&=U_{n}e^{\triangle t C}+\triangle t\phi_1(\triangle t C)F_n+\triangle t\phi_2(\triangle t C)(F_n-F_{n-1})\\
	&=U_{n}+\triangle t\phi_1(\triangle t C)(F_n+CU_n)+\triangle t\phi_2(\triangle t C)(F_n-F_{n-1}).
\end{align*}

To derive even higher-order schemes, one can build up higher-order approximations of the integrand as shown in the equation (\ref{Eq03}), i.e., with a reminder term of $O(\triangle t^q)(q\geq 4)$. For example, 
one can approximate $F(U(t_n+\tau), t_n+\tau)$ by high-order Taylor expansion and substitute into the integral term in the equation (\ref{Eq03}), leading to a family of $\phi$ functions (similar to $\phi_1$)
\begin{equation}\label{Eq07}
	\phi_k(\triangle t C)=\frac{1}{\triangle t^k}\int_{0}^{\triangle t}e^{(\triangle t-\tau)C}\tau^{k-1}d\tau=\int_{0}^{1}e^{(1-\lambda)\triangle t C}\lambda^{k-1}d\lambda,
\end{equation}
which are bounded satisfying the following recursion relation
\begin{equation}\label{Eq08}
	\phi_{k+1}(z)=\frac{k\phi_k(z)-1}{z},\; k\geq 1,\; {\rm where} \;\phi_0(z)=e^{z}.
\end{equation}
As discussed in \cite{cox2002exponential}, the explicit multistep ETD schemes with arbitrary order have been derived based on a polynomial approximation of  $F(U(t_n+\tau),t_n+\tau)$, 

\begin{equation}\label{Eq08b}
	U_{n+1}=U_ne^{\triangle t C}+\triangle t \sum_{m=0}^{s-1}g_m\sum_{k=0}^{m}(-1)^k\icol{m\\k}F_{n-k},
\end{equation}
where 
\begin{equation*}
	g_m=(-1)^m \int_0^1 e^{(1-\lambda)\triangle t C}\icol{-\lambda \\ m}d\lambda, \quad \text{with}    \quad \icol{-\lambda \\ m}=\frac{(-\lambda)(-\lambda-1)\cdots (-\lambda-m+1)}{m!}.
\end{equation*}

\subsection{Exponential time differencing  with Runge-Kutta time stepping}
As mentioned in \cite{cox2002exponential},  since the multistep explicit ETD schemes require $s$ previous evaluations of the nonlinear term $F$ as depicted in (\ref{Eq08b}), they are sometimes inconvenient to use. 
By adopting Runge-Kutta(RK) type approach alternatively, this inconvenience can be avoided. In addition,   Runge-Kutta(RK) type approaches typically have the advantages of smaller error constants and larger stability regions than the multistep explicit ETD methods. For instance, a brief summary of ETDRK schemes up to the fourth order is listed in the following:
\begin{itemize}
	\item ETD2RK
	
\begin{align*}\label{ETDRK2}
		a_n&=U_n+\triangle t(F_n+CU_n)\phi_1(\triangle t C)\\
	U_{n+1}&=U_{n}+\triangle t(F_n+CU_n)\phi_1(\triangle t C)+\triangle t(F(a_n, t_n+\triangle t)-F_{n})\phi_2(\triangle t C)
\end{align*}
	
	\item ETD3RK

	\begin{align*}
			a_n&=U_n+\frac{\triangle t}{2}(F_n+CU_n)\phi_1(\triangle t C/2)\\
		    b_n&=U_n+\triangle t\phi_1(\triangle t C)(CU_n+2F(a_n, t_n+\triangle t/2)-F_n)\\
		U_{n+1}&=U_{n}+\triangle t  \phi_1(\triangle t C)(CU_n+F_n)+\triangle t \phi_2(\triangle t C)(-3F_n\\
		&+4F(a_n, t_n+\triangle t/2)-F(b_n, t_n+\triangle t))\\
		&+\triangle t \phi_3(\triangle t C)(2F_n-4F(a_n, t_n+\triangle t/2)+2F(b_n, t_n+\triangle t))
	\end{align*}
	\item ETD4RK
	\begin{align*}
	a_n&=U_n+\frac{\triangle t}{2}(F_n+CU_n)\phi_1(\triangle t C/2)\\
	b_n&=U_n+\frac{\triangle t}{2}\phi_1(\triangle t C/2)(F(a_n, t_n+\triangle t/2)+CU_n)\\
	c_n&=a_n+\frac{\triangle t}{2}\phi_1(\triangle t C/2)(2F(b_n, t_n+\triangle t/2)-F_n+Ca_n)\\
	U_{n+1}&=U_{n}+\triangle t \phi_1(\triangle t C)(CU_n+F_n)\\
	&+\triangle t \phi_2(\triangle t C)(-3F_n+2(F(a_n, t_n+\triangle t/2)+F(b_n, t_n+\triangle t/2))-F(c_n, t_n+\triangle t))\\
	&+\triangle t \phi_3(\triangle t C)(2F_n-2(F(a_n, t_n+\triangle t/2)+F(b_n, t_n+\triangle t/2))+2F(c_n, t_n+\triangle t)).
	\end{align*}
\end{itemize}

\subsection{Computing a linear combination of $\phi$-functions}

The implementation of ETD schemes only requires computing the action of matrix function $\phi_{k}(A)$ on vectors $v_k$. For the evaluation of linear combinations of $\phi$-functions acting on sets of vectors $v_0$, $v_1$, $\cdots$, $v_p$, 
\begin{equation}\label{Eq09}
	\phi_0(A)v_0+\phi_1(A)v_1+ \cdots+\phi_{p}(A)v_p,
\end{equation}
it is crucial within calculations of all ETD schemes described above.

A few state-of-the-art algorithms to efficiently evaluate linear combinations of matrix function $\phi_{k}(A)$ on vectors $v_k$ are presented in \cite{luan2019further, gaudreault2016efficient, niesen2012algorithm, hochbruck1997krylov, al2011computing, caliari2016leja}. In this paper, we employ the algorithm $phipm_{-}simul_{-}iom2$ in \cite{luan2019further} to  evaluate the linear combination as in the equation (\ref{Eq09}), which typically consists of $(p+1)$ matrix-vector 
multiplications that can be carried out in a lower dimensional Krylov subspace. As a result, computational cost can be reduced significantly.

The linear combination (\ref{Eq09}) is actually equivalent to the solution of the following ODE 
\begin{equation*}\label{Eq10}
	y'(t)=Ay(t)+v_1+tv_2+\cdots+\frac{t^{p-1}}{(p-1)!}v_p, \; y(0)=v_0.
\end{equation*}

i.e., 
\begin{align*}
	y(1)&=e^Av_0 +\int_{0}^1 e^{(1-\lambda)A}d\lambda v_1+ \int_{0}^1 e^{(1-\lambda)A} \lambda d\lambda v_2+\cdots +\int_{0}^1 e^{(1-\lambda)A} \lambda^{p-1}d\lambda v_p, \\
	&=\phi_0(A)v_0 +\phi_1(A) v_1+ \phi_2(A) v_2+\cdots + \phi_p(A) v_p. 
\end{align*}

\section{Numerical Experiments}
In this section, the accuracy, stability and efficiency of the developed numerical methods as described above  are investigated through various testing examples. First we start with the accuracy test of the embedded boundary method by solving a Poisson equation with a virus-shaped geometry.  Next numerical experiments on the reaction-diffusion equation in irregular domains are presented to exhibit the accuracy, stability and efficiency of ETD schemes compared to other methods. Finally, we present the performance of the second-order ETD2 scheme combined with the level set method to solve a free boundary problem. 

\subsection{Convergence study for the embedded boundary method}
In this example, we consider solving a Poisson equation  $\nabla (\beta \nabla u )=f$ in an irregular domain $\Omega$ determined by the  boundary interface which is parameterized by 
\[
\begin{cases}
	x(\theta)=(0.6+0.1\sin(12\theta))\cos(\theta),\\
	y(\theta)=(0.6+0.05\sin(12\theta))\sin(\theta),
\end{cases}
\]
with $\theta \in [0,2\pi)$. The exact solution on $\Omega$ for this case is $u=e^x(x^2\sin(y)+y^2)$ and $\beta=2+\sin(xy)$.
A $(N+1)\times (N+1)$ uniform mesh partitioning $[-1,1]\times[-1,1]$ is used.

\begin{figure}[H]
	\centering
	\includegraphics[width=0.495\textwidth,trim={0.2cm 0.2cm 1.75cm 0.2cm},clip]{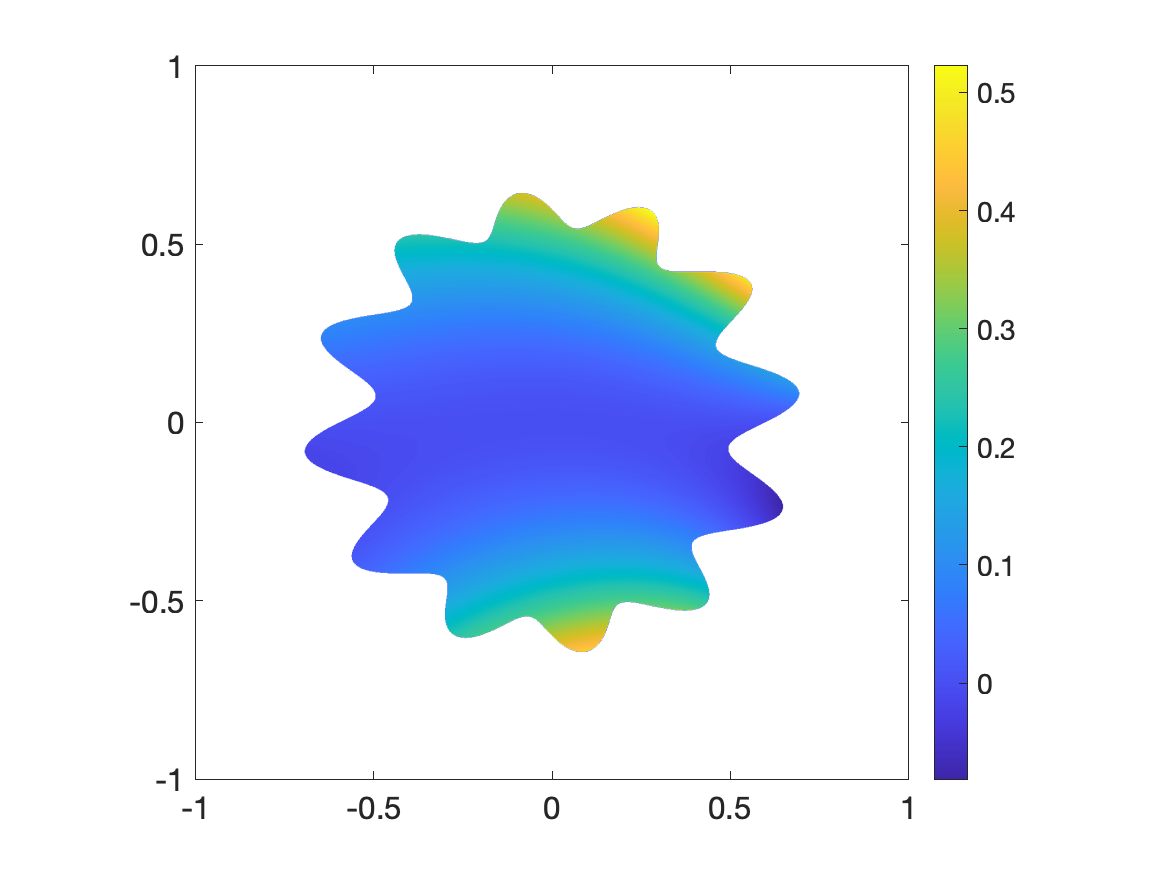}
	\includegraphics[width=0.495\textwidth,trim={0.2cm 0.2cm 1.75cm 0.2cm},clip]{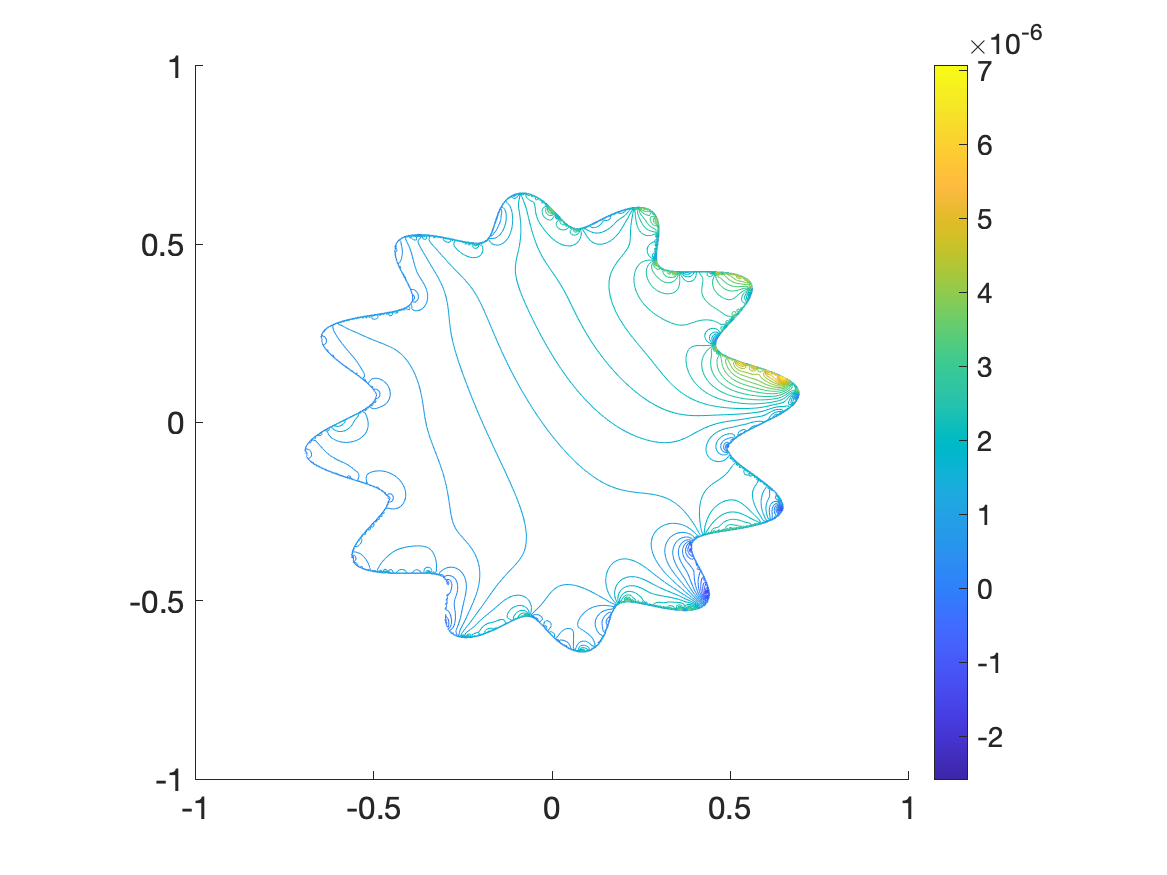}
	\caption{Numerical solution and error of the Poisson equation on a virus-shape geometry with grid points $1280\times1280$. Left: numerical solution. Right: numerical error.\label{fig:virus_sol}}
\end{figure}

\begin{figure}[h!]
	\centering
	\includegraphics[width=0.47\textwidth,trim={0.2cm 0.2cm 1.75cm 0.2cm},clip]{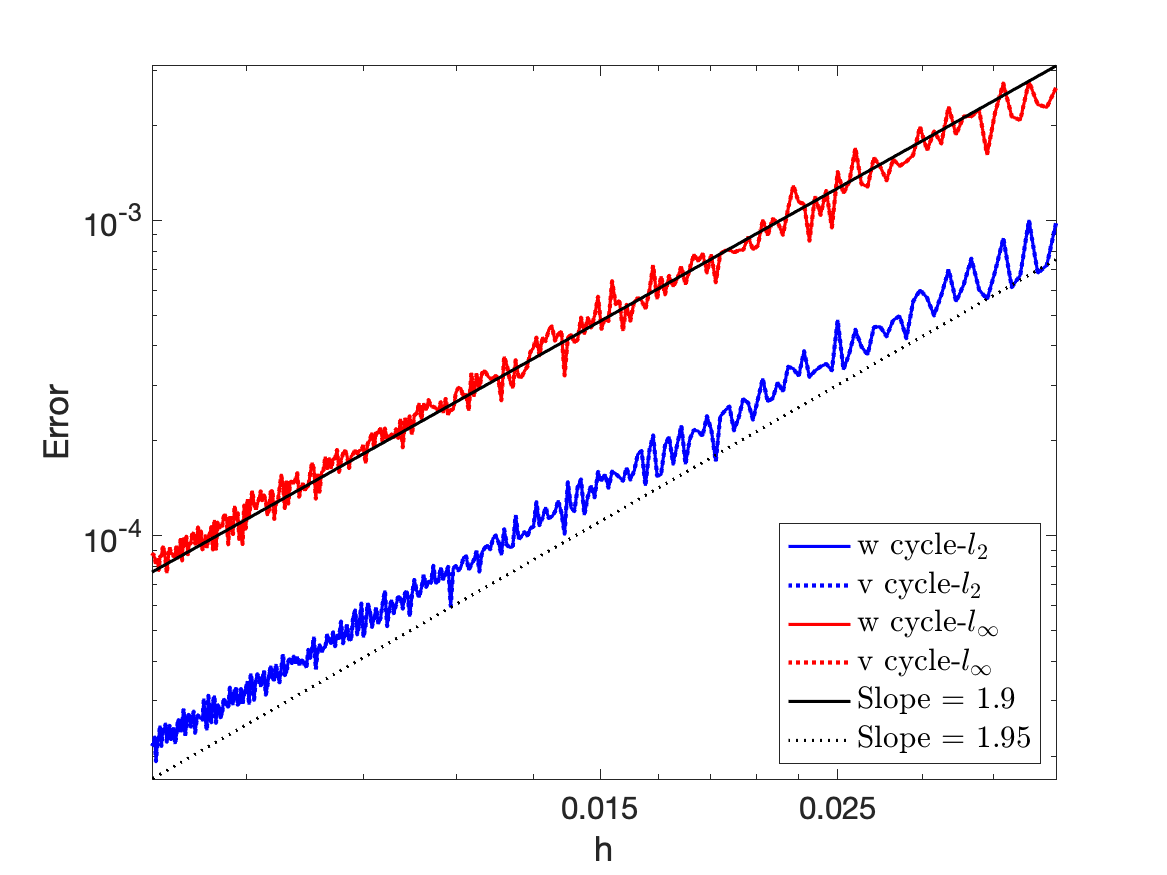}
	\includegraphics[width=0.47\textwidth,trim={0.2cm 0.2cm 1.75cm 0.2cm},clip]{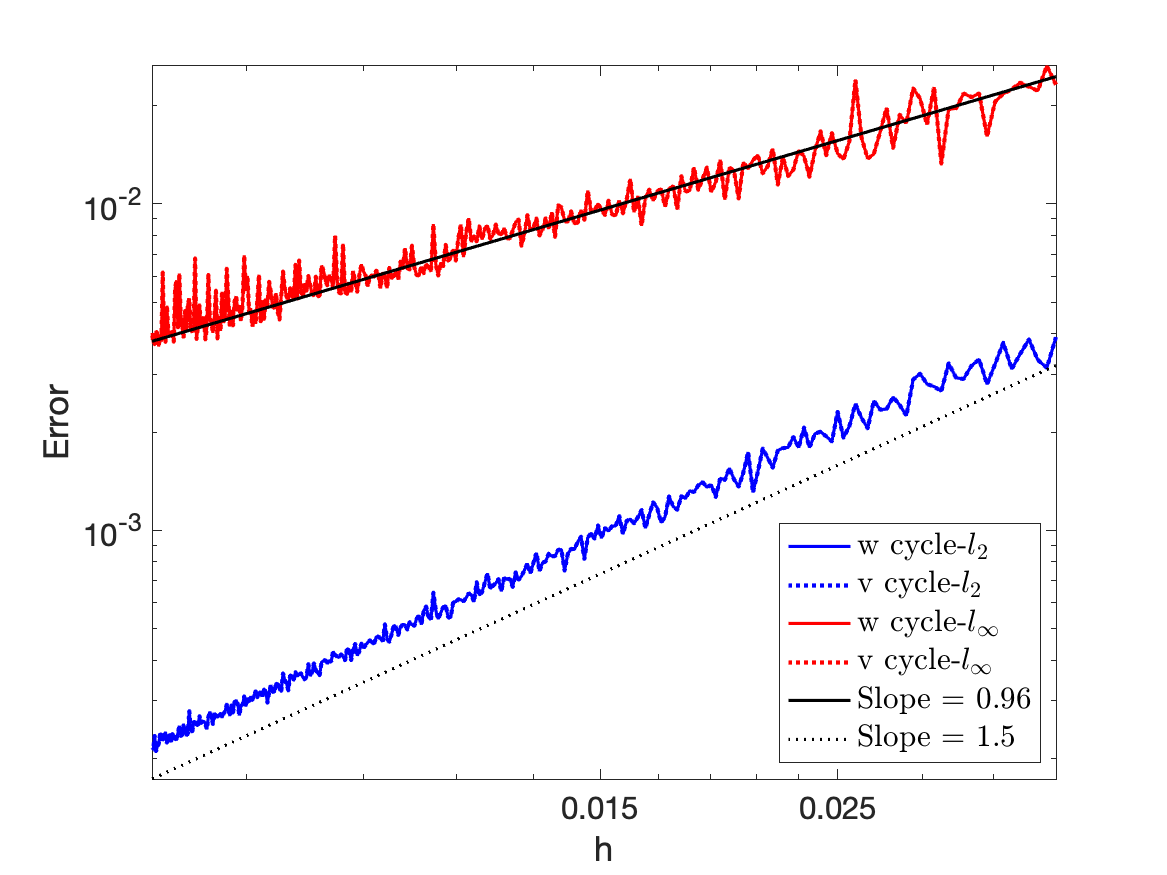}
	
	\caption{ Error analysis and convergence analysis of the solution with virus-shaped geometry using algebraic multigrid ``W" shape and ``V" shape.  Left: The numerical error of the solution. Right: The numerical error of the gradient of the solution. 	\label{fig:virus_error}}
\end{figure}

The numerical solution with $N=1280$ is presented in \fig{\ref{fig:virus_sol}}. Sweeping from $N=50$ to $N=350$, the second order accuracy for approximation of the solution can be observed 
with the developed embedded boundary method both in $l_2$ norm and $l_{\infty}$ norm (see \fig{\ref{fig:virus_error}} (left)). Furthermore, we can also observe $O(h^{1.5})$ accuracy in $l_2$ norm and $O(h^{0.96})$ accuracy in $l_{\infty}$ norm for approximation of the gradients of the solution (see \fig{\ref{fig:virus_error}} (right)). 

\subsection{Numerical tests for ETD with reaction-diffusion systems}
In this section, we incorporate the embedded boundary method for spatial discretization with ETD schemes to solve systems of reaction-diffusion equations in irregular domains. 
Numerical experiments are performed to demonstrate the accuracy, stability and efficiency of ETD schemes. Without loss of generality and for the convenience of better comparison, 
the test example is selected with analytical solutions given. 
Specifically,  we consider the following example with a peanut-shaped and non-convex geometry which is determined by the level-set function. 
\begin{equation*}
	\rho(x,y)=0.5-e^{-20(x^2+(y-0.25)^2)}-e^{-20(x^2+(y+0.25)^2)}.
\end{equation*}
Here we solve the reaction-diffusion equation $u_t=\nabla (\beta \nabla u )+f$ defined inside the domain $\rho(x,y)<0$, where $\beta =0.25-x^2-y^2$. The source term $f$ is computed by assuming a exact solution $u(x,y,t)=e^{-t}(x^2+y^2-0.25)$. The computational domain is $[-1,1] \times [-1, 1]$. In order to make a fair comparison between different methods, for this example, we will mainly focus on the following second-order numerical schemes: the second-order Crank-Nicolson, ETD2, and ETD2RK. 
 
\begin{figure}[h]
	\centering
	\includegraphics[width=0.495\textwidth,trim={0.2cm 0.2cm 1.75cm 0.2cm},clip]{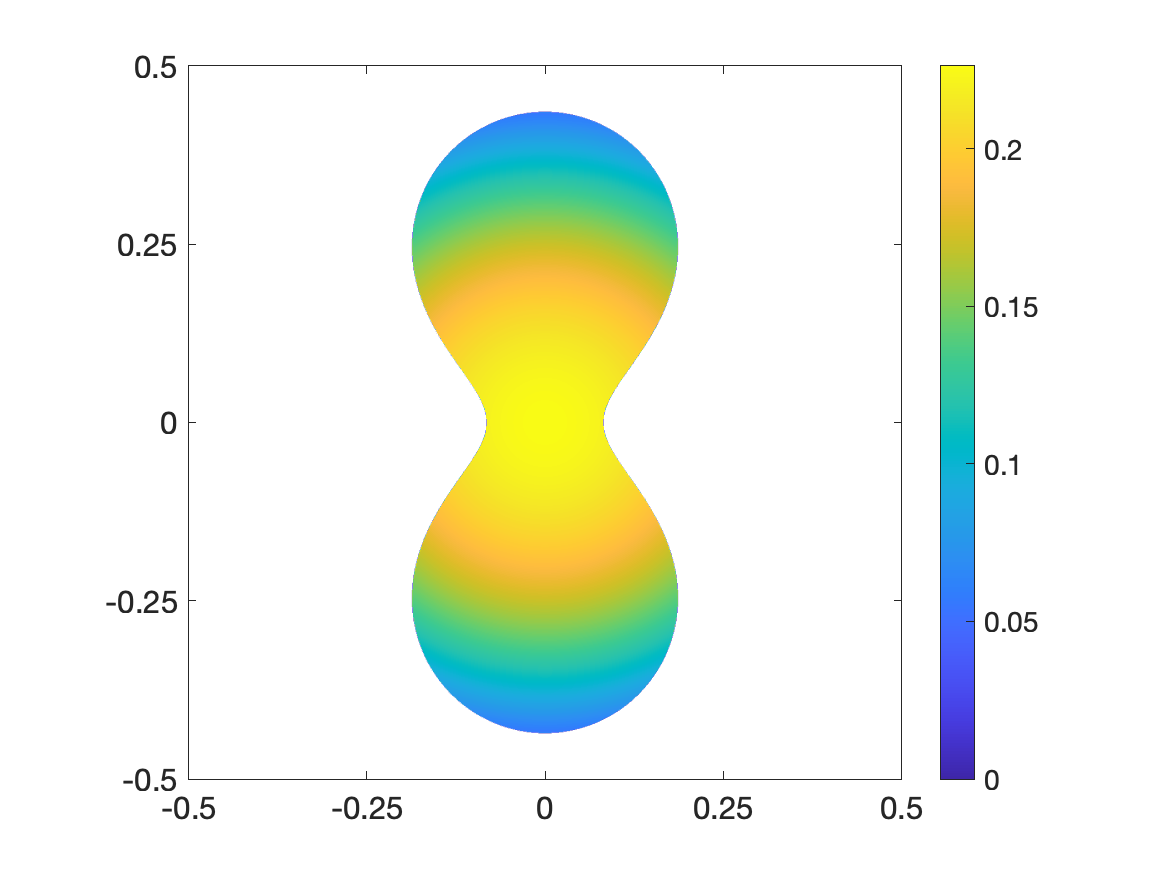}
	\includegraphics[width=0.495\textwidth,trim={0.2cm 0.2cm 1.75cm 0.2cm},clip]{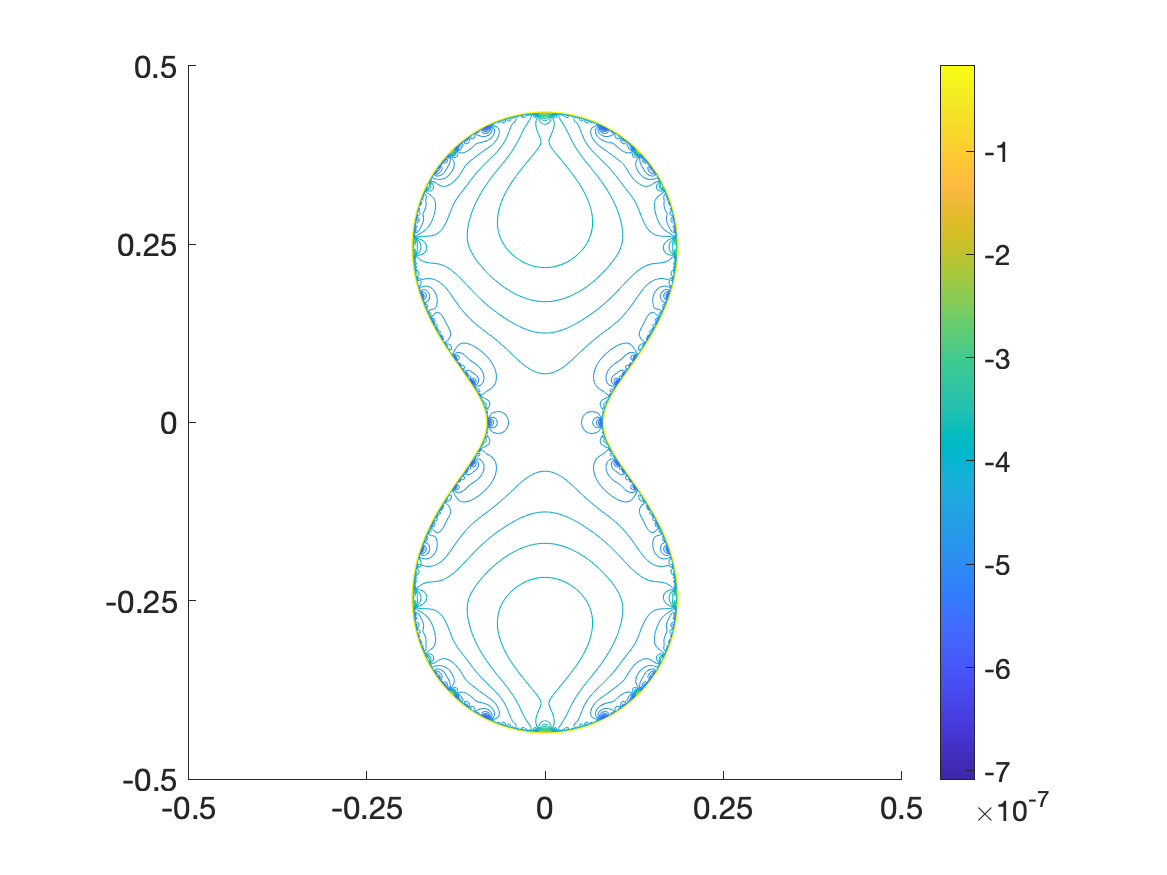}
	\caption{Numerical solution and error  of the reaction-diffusion equation at $t=0.5$ with a peanut-shape geometry with grid points $1280\times1280$ using ETD2. Left: numerical solution. Right: numerical error. \label{fig:peanut_sol}}
\end{figure}

\subsubsection{Accuracy test} 

Numerical errors and corresponding convergence rates of the second-order Crank-Nicolson method, ETD2 and ETD2RK at time $t=0.1$ with five different spatial and temporal resolutions are reported in Table \ref{table:1}, where the time step $dt$ are taken as equal to the grid size $h$. As expected, we can clearly see  a second-order accuracy for all three schemes. For illustration, the numerical solution and numerical error of the reaction-diffusion equation with ETD2 at final time $t=0.5$ with grid points 1280 $\times$ 1280 are presented in Figure \ref{fig:peanut_sol}.

\begin{table}[H]
\centering
\setlength{\abovecaptionskip}{-6pt}
\caption{Numerical errors in $l_{\infty}$ norm  and $l_2$ norm and corresponding convergence rates for Crank-Nicolson, ETD2, and ETD2RK schemes at $t=0.1$. }	\label{table:1}
\begin{center}
\begin{tabular}{ |c||c|c|c|c|  }
\hline
	\multicolumn{5}{|c|}{Convergence test for Crank-Nicolson  }\\
\hline
$N_x\times N_y \times N_t$ & $l_{\infty}$-Error&Order&$ l_2 $-Error&Order\\
\hline
81$\times$81$\times $4&8.041$\times10^{-4}$   &- &2.296$\times10^{-4}$&- \\
\hline
161$\times$161$\times $8&1.748$\times10^{-4}$   & 2.201&4.436$\times10^{-5}$&2.372\\
\hline
321$\times$321$\times $16&5.314$\times10^{-5}$  &1.718&1.059$\times10^{-5}$ &2.067\\
\hline
641$\times$641$\times $32&1.859$\times10^{-5}$    & 1.515&2.751$\times10^{-6}$&1.945\\
\hline
1281$\times$1281$\times $64&4.333$\times10^{-6}$     & 2.101&6.254$\times10^{-7}$ &2.137\\
\hline
\hline
\multicolumn{5}{|c|}{Convergence test for ETD2 }\\
\hline
$N_x\times N_y \times N_t$ & $l_{\infty}$-Error&Order&$ l_2 $-Error&Order\\
\hline
81$\times$81$\times $4&8.234$\times10^{-4}$   &- &2.580$\times10^{-4}$&- \\
\hline
161$\times$161$\times $8&1.887$\times10^{-4}$   & 2.126&5.542$\times10^{-5}$&2.219\\
\hline
321$\times$321$\times $16&5.759$\times10^{-5}$  &1.712&1.420$\times10^{-5}$ &1.964\\
\hline
641$\times$641$\times $32&1.956$\times10^{-5}$    & 1.558&3.678$\times10^{-6}$&1.949\\
\hline
1281$\times$1281$\times $64&4.583$\times10^{-6}$     & 2.094&8.780$\times10^{-7}$ &2.067\\
\hline
\hline
\multicolumn{5}{|c|}{Convergence test for ETD2RK }\\
\hline
$N_x\times N_y \times N_t$ & $l_{\infty}$-Error&Order&$ l_2 $-Error&Order\\
\hline
81$\times$81$\times $4&7.851$\times10^{-4}$   &- &2.276$\times10^{-4}$&- \\
\hline
161$\times$161$\times $8&1.756$\times10^{-4}$   & 2.161&4.465$\times10^{-5}$&2.350\\
\hline
321$\times$321$\times $16&5.330$\times10^{-5}$  &1.720&1.071$\times10^{-5}$ &2.060\\
\hline
641$\times$641$\times $32&1.863$\times10^{-5}$    & 1.517&2.784$\times10^{-6}$&1.944\\
\hline
1281$\times$1281$\times $64&4.343$\times10^{-6}$     & 2.101&6.346$\times10^{-7}$ &2.133\\
\hline
\end{tabular}
\end{center}
\end{table}

\subsubsection{Stability test}
In this section, we test the stability properties of four numerical schemes: the standard explicit Runge-Kutta, ETD2, ETD2 Runge-Kutta, and Crank-Nicolson for solving the example of reaction-diffusion system with the peanut-shaped geometry as previously mentioned. 
We set the final time $t_{end}=0.2$ and a uniform grid size $h=0.004$ for all the simulations. The errors are measured in $l_2$ norm between numerical solutions and the exact solutions by varying time steps. 

\begin{figure}[H]
	\centering
	\includegraphics[width=0.65\textwidth,trim={0.1cm 0.1cm 1.33cm 0.cm},clip]{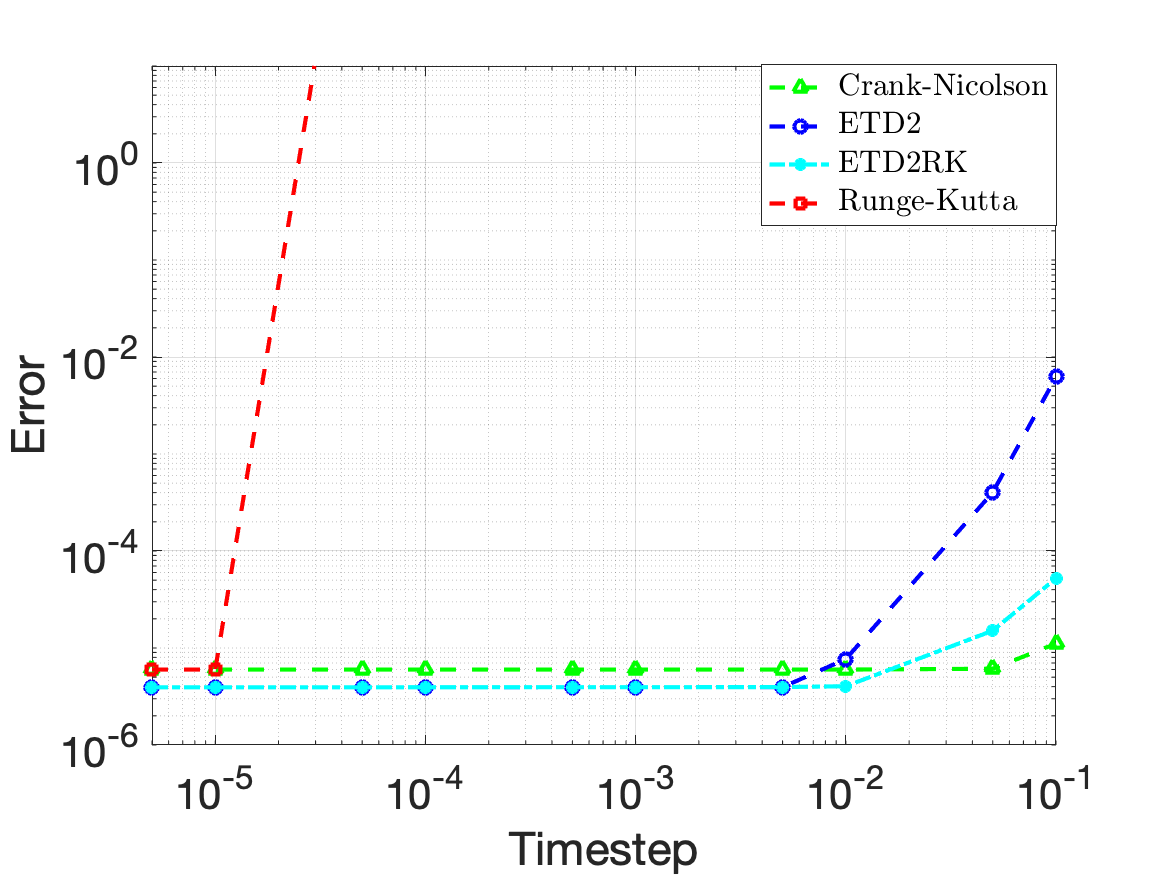}
	\caption{Errors between numerical solutions and the exact solution using different numerical schemes by varying time step sizes with a fixed uniform grid size $h=0.004$.}\label{stabilitytest}
\end{figure}

From Figure \ref{stabilitytest}, it can be observed that the explicit Runge-Kutta scheme blows up with time step size greater than $1\times 10^{-5}$ as expected. In contrast, all other three numerical schemes: Crank-Nicolson, ETD2 and ETD2RK 
exhibit very excellent stability conditions, which allow for very large time step size till $dt=0.1$.

\subsubsection{Efficiency test}
In Table \ref{table:4},  we compare the efficiency performance of three schemes with nice stability conditions: Crank-Nicolson, ETD2, and ETD2RK for a system with more refined grids. Here we choose the grid size $h=0.02$, $h=0.01$ and $h=0.005$, respectively. The time step size is  $dt=10^{-4}$ for all simulations. Note that when solving the linear system from the Crank-Nicolson method, we use the fast linear solver by the conjugate gradient method with an incomplete Cholesky preconditioner. 
However, solving a large size of linear system in each time step is still very costly. 
By avoiding solving large linear systems, it can be clearly observed that ETD2 is 2-3 times faster than Crank-Nicolson, in which we adopt the adaptive Krylov space to compute the multiplication of matrix and vectors.  
Since ETD2RK is a two-stage numerical algorithm, its efficiency is slightly better or comparable to Crank-Nicolson for this example.

\begin{table}[H]
	\centering
	\setlength{\abovecaptionskip}{-6pt}
	\caption{Efficiency test for a reaction-diffusion system by varying the number of grid points with the same time step size (Unit: seconds).}
	\label{table:4}
	\begin{center}
		\begin{tabular}{ lccc  }
			\hline
	         $n_x \times n_y$ & 1001 $\times$ 1001& 2001 $\times$ 2001 & 3001 $\times$ 3001\\
			\hline
			Crank-Nicolson &304.23& 2440.19&8183.33\\
			ETD2 &136.67&957.34&3176.54\\
			ETDRK2&275.89&2057.18&7325.82\\
			\hline
		\end{tabular}
	\end{center}

\end{table}

\subsection{Numerical tests of the free boundary problem}
For the free boundary problem as described in (\ref{eq:free}), it is very challenging and crucial to accurately and efficiently handle the reaction-diffusion equation with the changing domain for each time step. In this paper, we integrate the level set method to track the evolution of the moving boundary, and the ETD2 schemes with the discretization technique to solve the reaction-diffusion equation in each time step. As mentioned above, ETD2 schemes exhibit very nice stability conditions by allowing for large time step size, and it is also much faster than other schemes with the similar stability conditions like Crank-Nicolson and ETD2RK. Here we briefly introduce the numerical algorithm for solving a diffusive logistic model for the population of the invasive species $u(\mathbf{x},t)$  with free boundaries as follows,

\[
\begin{cases}
u_t=D\Delta u+u(a-bu), \quad \quad \quad \quad \quad \quad \quad \quad \quad t>0, \quad \mathbf{x}\in \Omega(t),\\
u(\mathbf{x},t)=0, \quad \quad \quad \quad \quad \quad \quad \quad \quad \quad \quad \quad \quad \quad t>0, \quad \mathbf{x}\in \Omega^c(t),\\
\vec{v}(\mathbf{x},t)=\mu |\nabla u(\mathbf{x},t)| \mathbf{n}(\mathbf{x})=-\mu \nabla u(\mathbf{x},t),~t>0, \quad \mathbf{x}\in \partial \Omega(t),\\

u(\mathbf{x},0)\in C^2(\Omega(0)), ~u(\mathbf{x},0)>0~\text{in}~\Omega(0),~ u(\mathbf{x},0)=0~\text{on}~\partial \Omega(0).
\end{cases}
\]
we introduce a level set function $\rho(x,y,t)$, such that $\rho(x,y,t)=0$ on the boundary, $\rho(x,y,t)<0 $ in $\Omega(t)$, and $\rho(x,y,t)>0 $ outside of $\Omega(t)$. $\rho(x,y,t)$ is initialized as a signed distance function to the initial boundary $\partial \Omega(0)$. 

{\bf Algorithm.}
\begin{compactenum}
	\item[Step 0.]
	Input all the parameters. 
	Set the computational box $ (-L, L)^2 \subset \R^2$ and cover it with a uniform
	finite-difference grid with grid sizes $h$. 
	Discretize the time interval $[0, T]$ of interest with time step $dt.$ 
	Initialize the level-set function $\rho^0$ and  $u^0$. 
	Set $m = 0.$ 
	
	\item[Step 1.]

	Extend the normal velocity $u(\mathbf{x},t_m)$ from the interface to the entire computational box.
	Discretize  the level set advection equation $	\rho_t+\vec{v} \cdot \nabla \rho=0$ with a HJ-WENO scheme. Solve it to get the updated level-set function $\rho^{m+1}$. 
	Reinitialize the level-set function and still denote it by $\rho^{m+1}$. 
	
	\item[Step 2.]
	Extend $u^m$  to new unknowns  overlapping with $\Omega_{m+1}$ defined by $\rho^{m+1}$ by a quadratic extrapolation in the normal direction by following \cite{aslam2004partial}.  Solve reaction diffusion equation with ETD schemes in the irregular domain $\Omega_{m+1}$ to obtain $u^{m+1}$.

	\item[Step 3.]

	Set $m:=m+1$. Repeat Steps 1--2 until the final simulation time is reached.  
	
\end{compactenum}

To show the potential application of the developed ETD2 for the Stefan-type free boundary problems, we consider one example from  \cite{liu2019numerical} with the following initial setups: $(D,\mu, a,b)=(1.5,1,1,1)$, with an initial square domain of length 0.5 centered at $(0,0)$. The initial level set function is
\begin{equation*}
	\rho(x,y,0)=-\text{min}(0.5-|x|).
\end{equation*}
and the initial function is 
\begin{equation*}
	u(x,y.0) = 
	\left\{
	\begin{array}{ll}
		20(0.5-x)^2(0.5+x)^2, & \text{if } (x,y)\in \Omega(0),\\
		0, & \text{if }(x,y)\in \Omega^c(0).
	\end{array}
	\right.
\end{equation*}

\begin{figure}[H]
	\centering
	\includegraphics[width=0.24\textwidth,trim={0.2cm 0.2cm 1.35cm 0.2cm},clip]{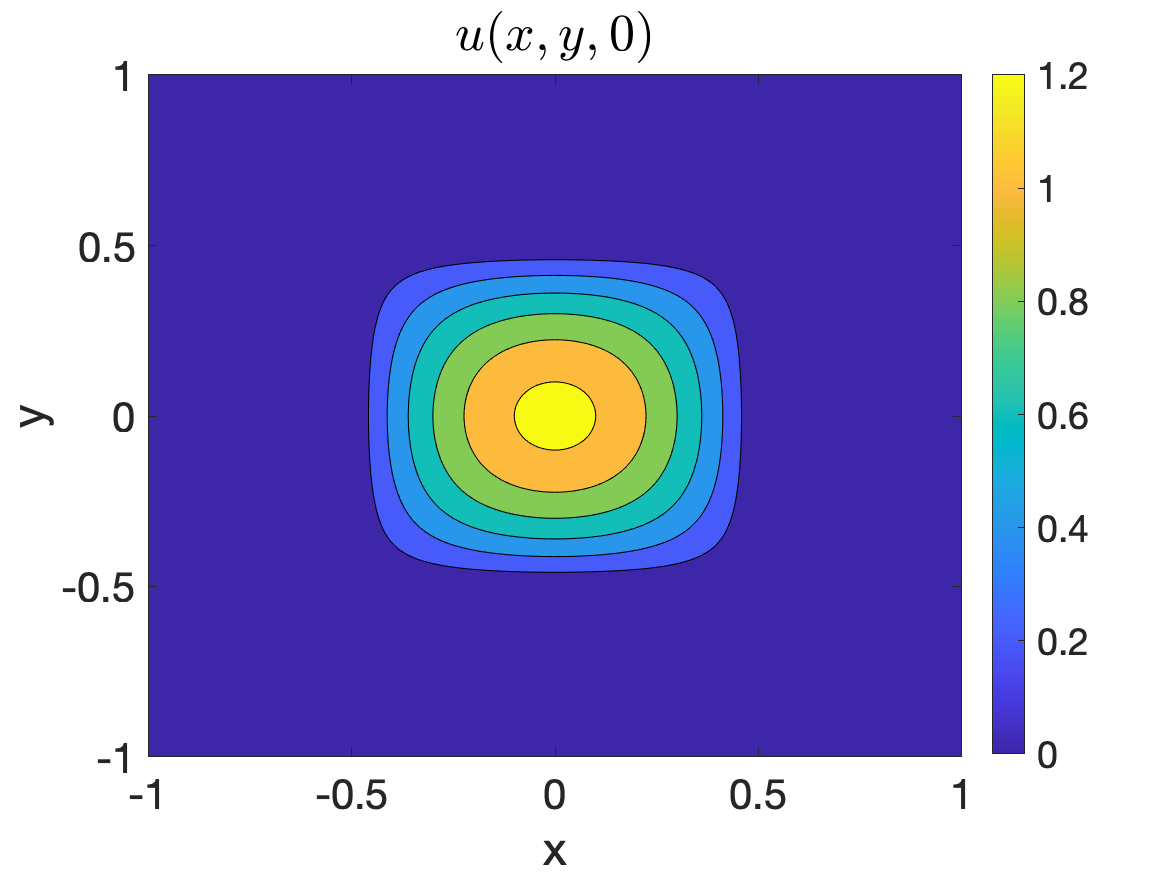}
	\includegraphics[width=0.24\textwidth,trim={0.2cm 0.2cm 1.35cm 0.2cm},clip]{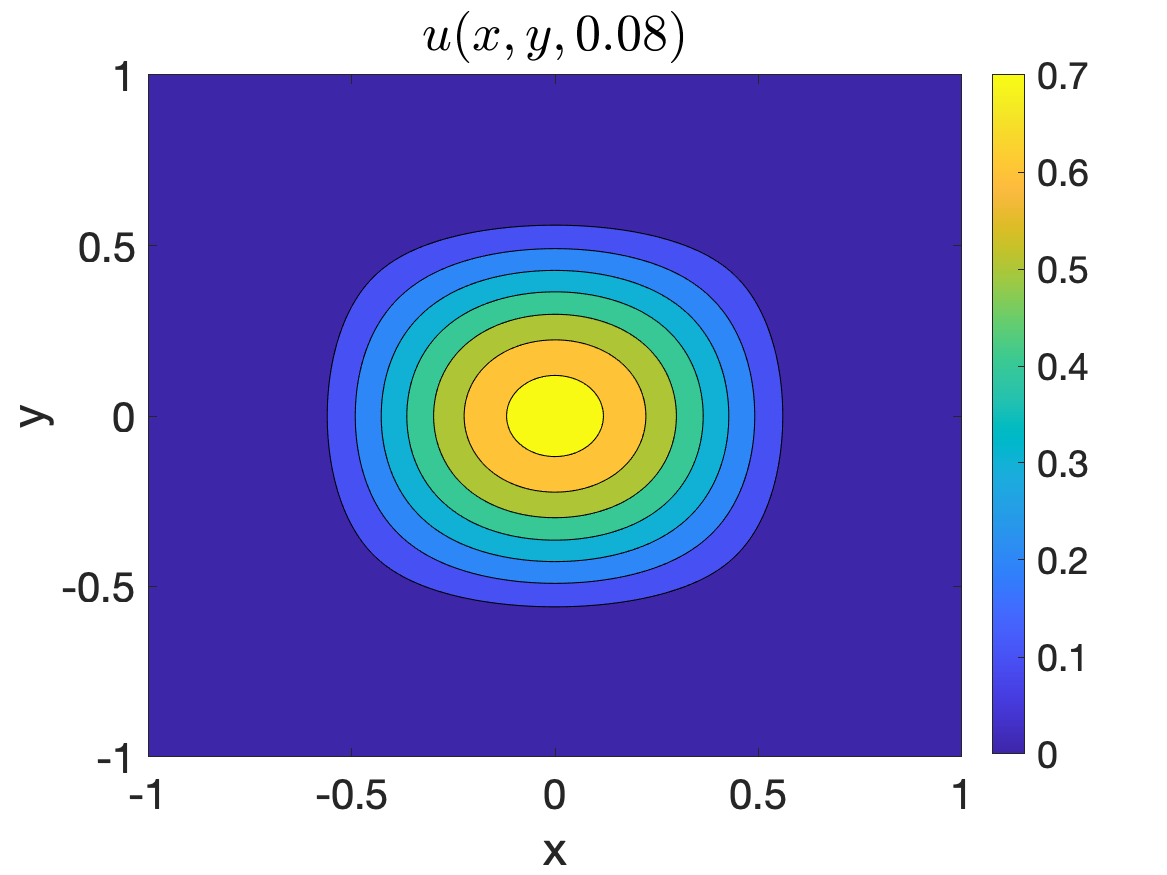}
	\includegraphics[width=0.24\textwidth,trim={0.2cm 0.2cm 1.0cm 0.2cm},clip]{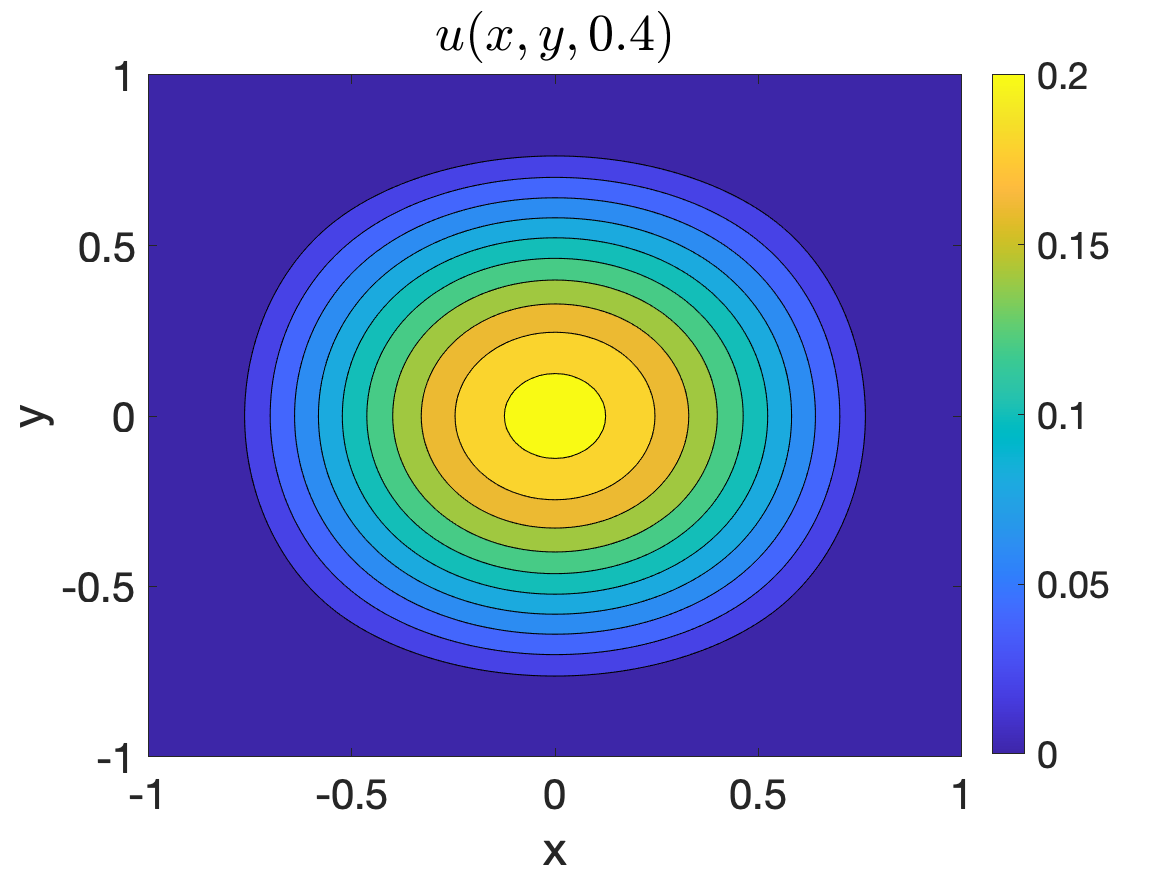}
	\includegraphics[width=0.24\textwidth,trim={0.2cm 0.2cm 0.65cm 0.2cm},clip]{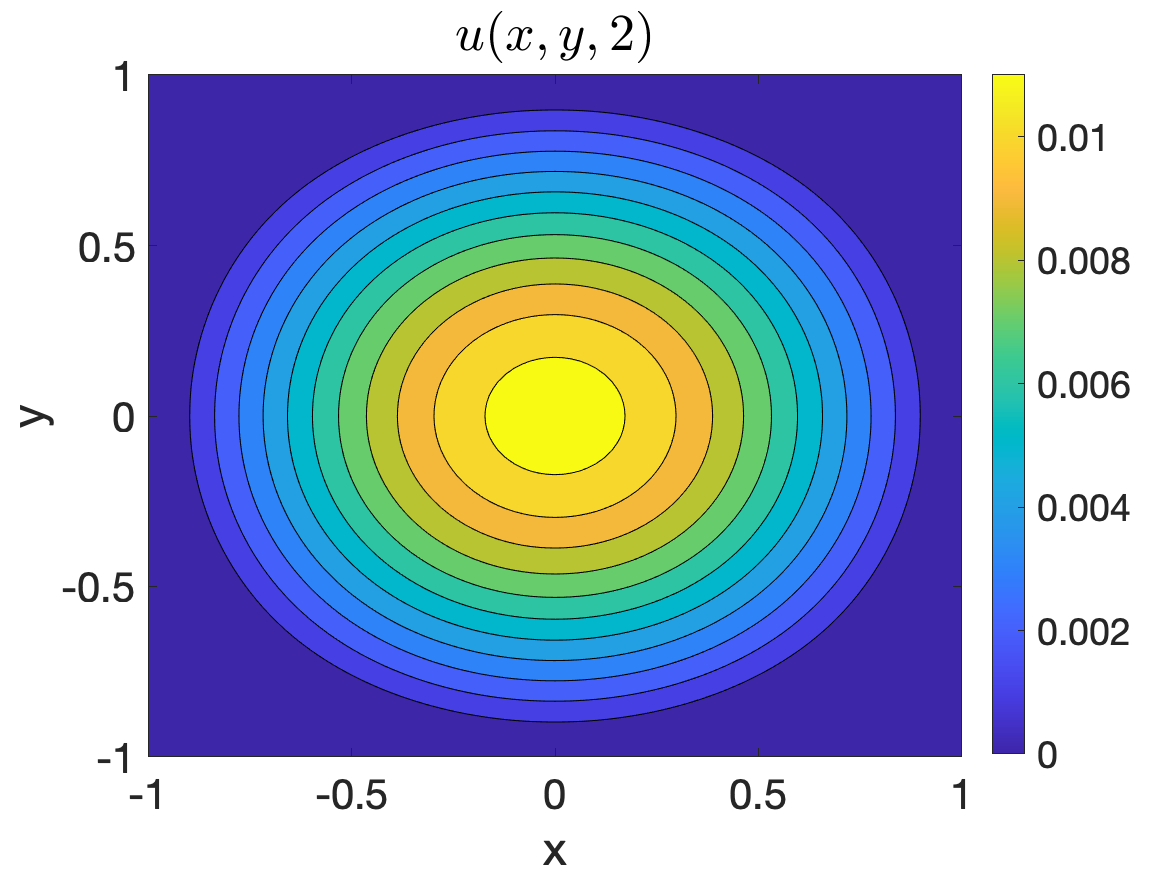}
	\includegraphics[width=0.24\textwidth,trim={0.2cm 0.2cm 1.75cm 0.2cm},clip]{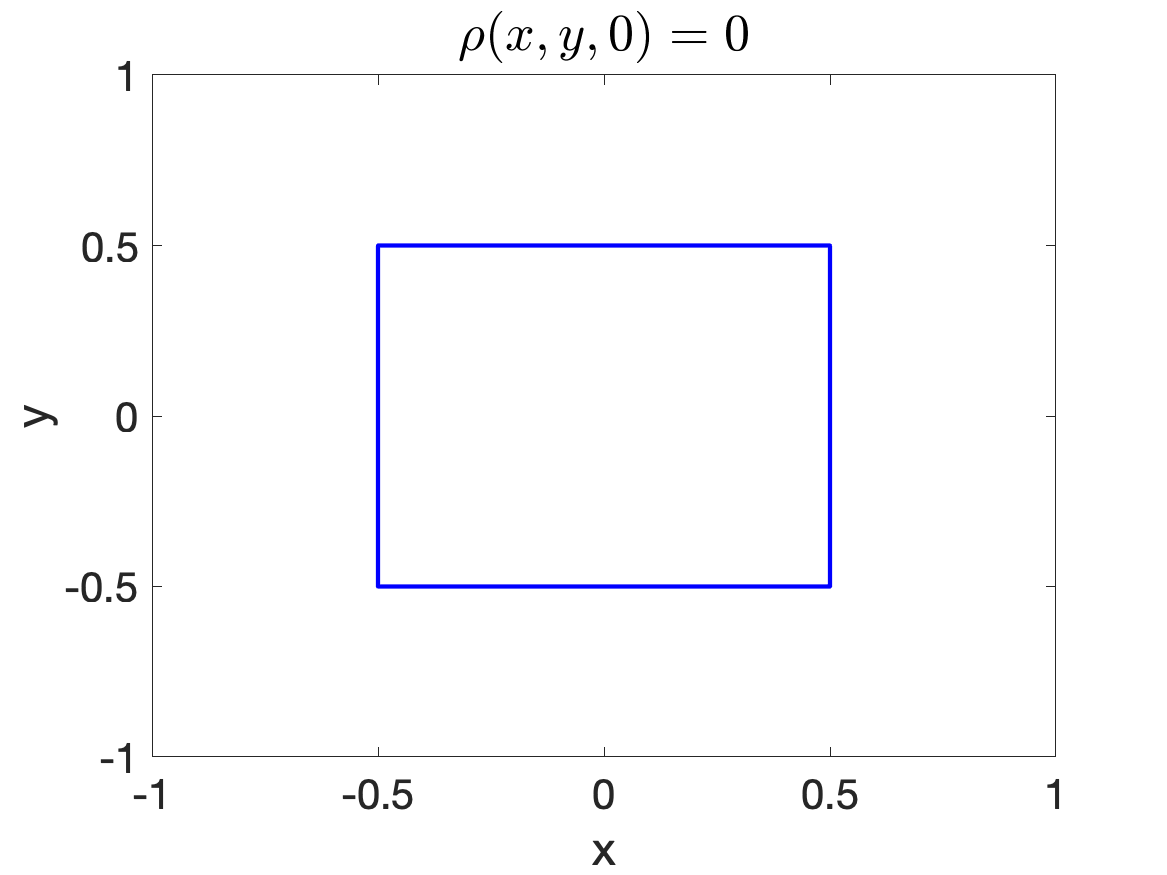}
	\includegraphics[width=0.24\textwidth,trim={0.2cm 0.2cm 1.75cm 0.2cm},clip]{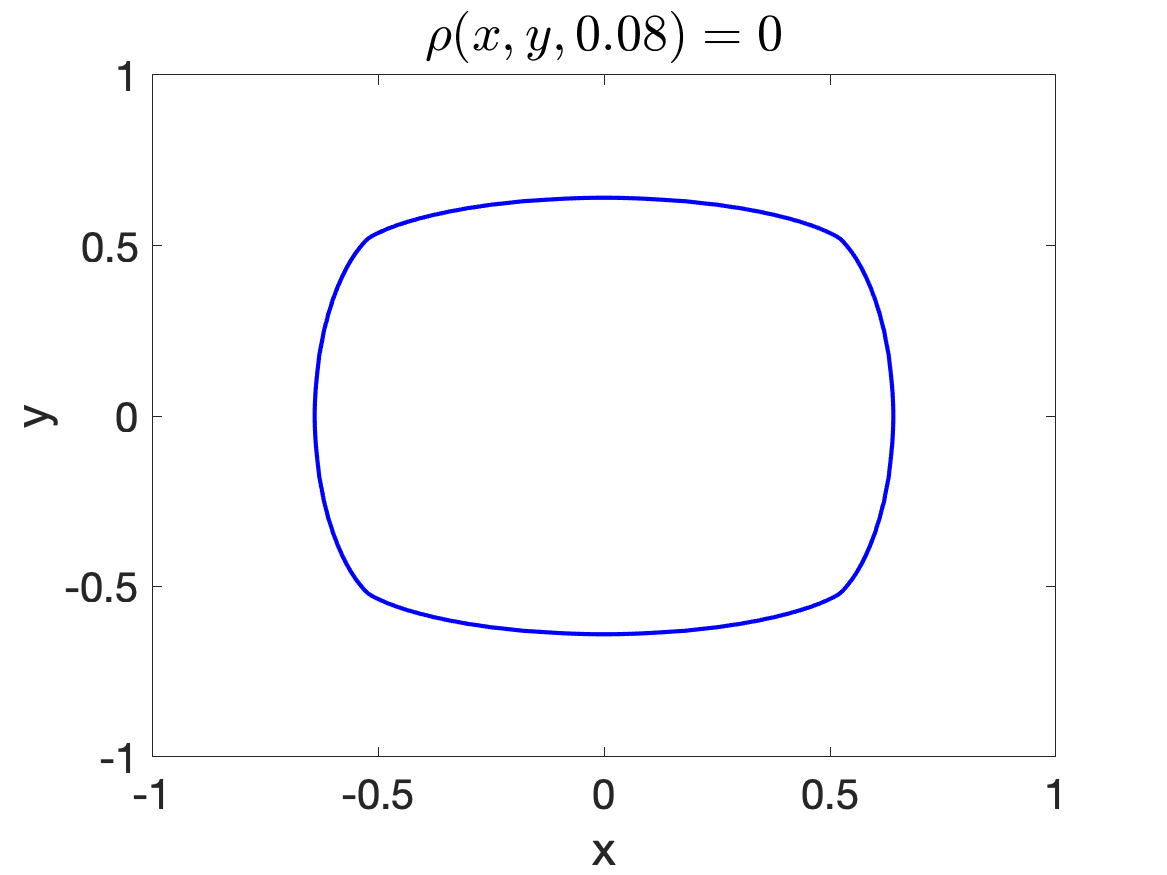}
	\includegraphics[width=0.24\textwidth,trim={0.2cm 0.2cm 1.75cm 0.2cm},clip]{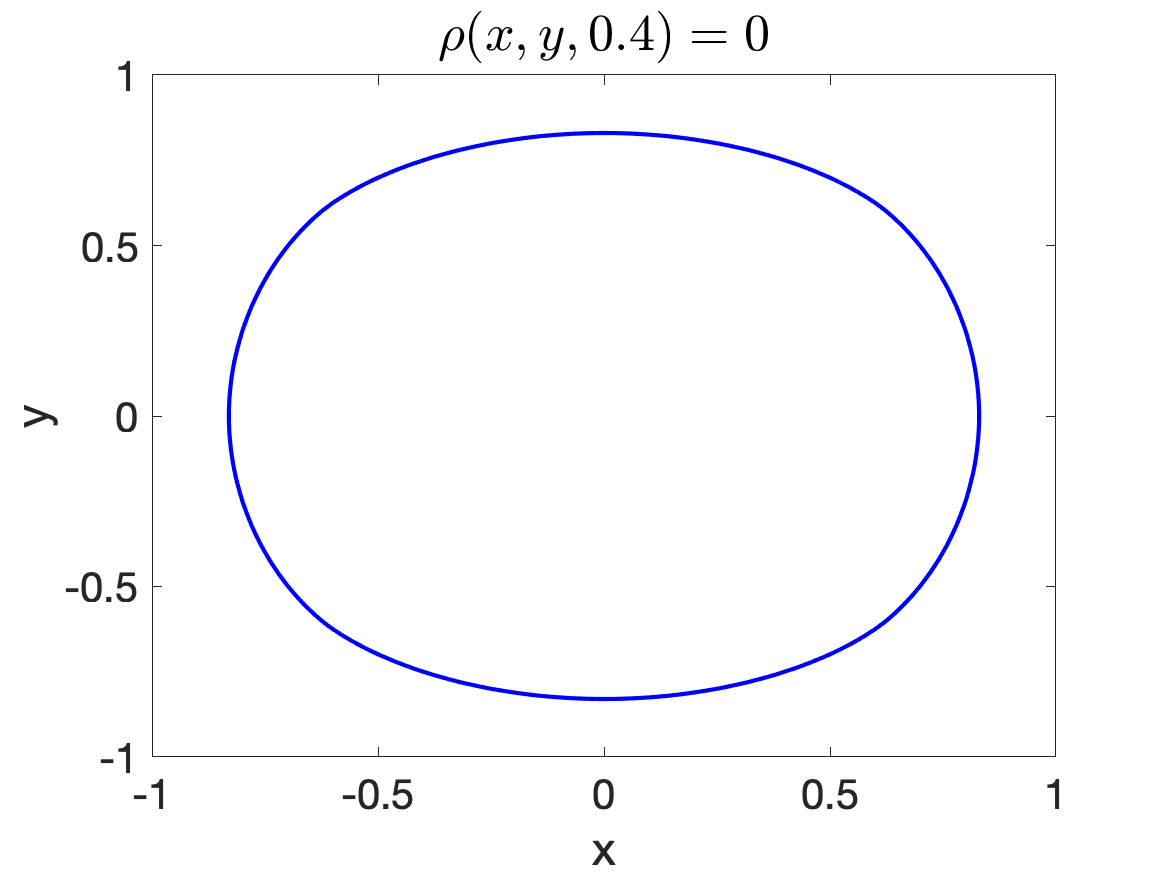}
	\includegraphics[width=0.24\textwidth,trim={0.2cm 0.2cm 1.75cm 0.2cm},clip]{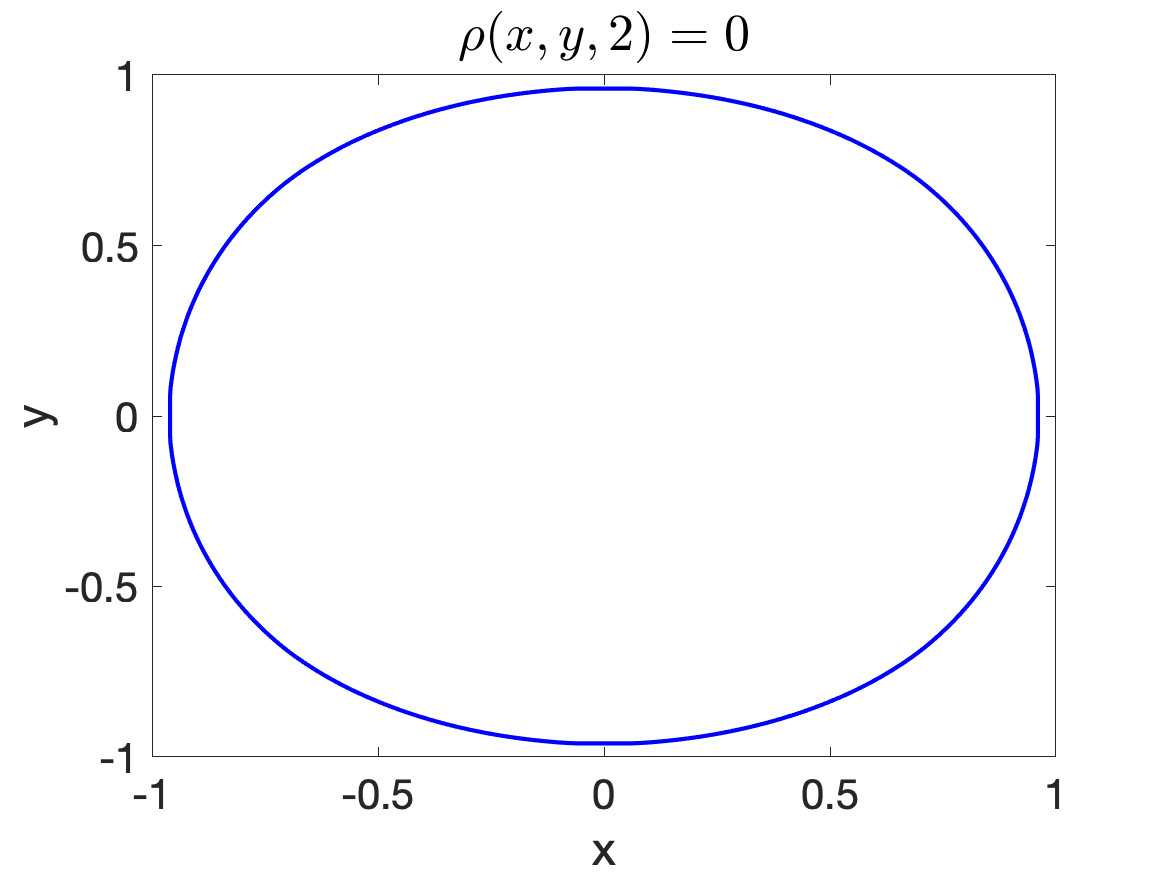}
	
	\caption{Evolution of $u(x,y,t)$ and the moving boundary $\rho(x,y,t)=0$ with the initial domain $\Omega(0)$  a square in 2D.  \label{fig:square}}
\end{figure}

  Figure~\ref{fig:square}~shows the evolution of the spreading of species  $u(x,y,t)$ along with the moving boundary. For this example, the moving boundary will asymptotically evolve into circles, which correlates exactly with the theoretical asymptotic behavior described in  \cite{du2014regularity}.

\section{Conclusion}
In this paper, we have incorporated the embedded boundary method, ETD schemes with level set method to systematically study reaction-diffusion systems in irregular domains with free boundaries. 
To our best knowledge, it is the first work to integrate ETD scheme with the embedded boundary method for time-dependent PDEs as well as to combine ETD with level set method for solving free boundary problems in two dimensions.
Through numerical experiments, we first show the accuracy of the embedded boundary method for a Poisson equation with a virus-shaped geometry. Next we test the accuracy, stability, and efficiency of the ETD schemes along with other methods. 
In order to significantly reduce the computational cost, we have adopted the state-of-the-art algorithm: $phipm_{-}simul_{-}iom2$  \cite{luan2019further} to evaluate the linear combination of matrix-vector multiplications in ETD schemes using a lower dimensional adaptive Krylov subspace. In summary, ETD scheme is superior to other three selected schemes (RK, Crank-Nicolson and ETDRK) in terms of a combination of accuracy, stability and efficiency, especially in efficiency.
More importantly, the ETD2 scheme has been successfully employed to a reaction-diffusion system with free boundaries, which produces very promising results for the free boundary problems.

Our immediate next step is to further develop a fast and efficient numerical algorithm for other more complicated systems with free boundaries. Possible extensions include (but not limited to) combining the fast local level set method and the ETD schemes for free boundary problems, and the generalization of the developed methods to the Navier-Stokes equations and the Grad-Shafranov equations.

\section*{Statement of no conflict of interest}

{{On behalf of all authors, the corresponding author states that there is no conflict of interest.}}

\bibliographystyle{plain}
\bibliography{SLbibfile}

\end{document}